\documentclass{amsart}
\usepackage{amsmath}  

\def\eg{{\it e.g.\ }} \def\ie{{\it i.e.\ }}
\def\Real{{\mathbb R}}

\theoremstyle{plain}

\theoremstyle{remark}

\theoremstyle{definition}

\begin{document}
\pagestyle{myheadings}
\markright{Fabian Waleffe \hfill  September 14, 2004 \hfill}

\title{Remarks on a quasi-linear model of the Navier-Stokes equations}
\author{Fabian Waleffe}
\address{Department of Mathematics,
University of Wisconsin, Madison, WI 53706}
\email{waleffe@math.wisc.edu}
\date{\today}

\begin{abstract}
Dinaburg and Sinai recently proposed a quasi-linear model of the
Navier-Stokes equations. Their model assumes that nonlocal interactions
in Fourier space are dominant, contrary to the Kolmogorov 
turbulence phenomenology where local interactions prevail. 
Their equation corresponds to the linear evolution
of small scales on a background field with uniform gradient,
but the latter is defined as the linear superposition of all
the small scale gradients at the origin. This is not self-consistent.
\end{abstract}

\maketitle

Dinaburg and Sinai \cite{DS01} recently proposed a quasi-linear approximation
of the Navier-Stokes equations which they feel preserves the
basic character of the Navier-Stokes nonlinearity.
They prove existence and uniqueness of solutions to their model,
for special cases \cite{DS03}. 
Here, we show that their equation is identical to that governing
the linear evolution of small scales on an infinitely large 
scale flow with uniform gradient (eqn.\ (\ref{FTveq}) below).
This is a direct consequence of their assumption that nonlocal
interactions in Fourier space dominate.
That assumption is contrary to a large body of phenomenological,
experimental and numerical studies of the turbulent energy cascade 
 where local interactions are thought to dominate
(see \eg\cite{Frisch}).  On that basis alone, it seems
 unlikely that  the model would preserve the basic character 
of the Navier-Stokes nonlinearity, but the model also 
contains a basic inconsistency.  It defines the velocity gradient
 of the large scale field as
the net velocity gradient induced by all
 the small scales at the origin (eqn.\ (\ref{ADS}) below).
This is not self-consistent since small
scales do not uniformly distort larger scales.
A self-consistent mean field theory, where 
the large scale flow results from the nonlinear
interactions of the small scale fluctuations, instead of 
their linear superposition, is not possible in the nonlocal
limit.


Consider the Navier-Stokes equations for incompressible flow of a
viscous fluid in three-dimensional space $\Real^3$:
\begin{equation}
\begin{split}
\frac{\partial u}{\partial t} + (u\cdot \nabla ) u + \nabla p =&\;  \nu
\nabla^2 u,\\
\nabla \cdot u =&\; 0, \quad
\label{eqn:NSE}
\end{split}
\end{equation}
where $u=u(x,t)=(u_1(x,t),u_2(x,t),u_3(x,t))$ is the velocity vector
 field,  $p=p(x,t)$ is the kinematic pressure, $x=(x_1,x_2,x_3)$ is 
the position vector and $\nu>0$ is the kinematic viscosity.

Let $u=A\cdot x + v$ where $A=A(t)$ is a matrix in $\Real^3
\times \Real^3$ such that $dA/dt+ A\cdot A$ is symmetric 
and $\mbox{tr}(A)=0$ in order for $A\cdot x$ to be a
 solution of the incompressible Navier-Stokes
equations (\ref{eqn:NSE}). The base field $A\cdot x$ has
uniform gradient and $v=v(x,t)=(v_1(x,t), v_2(x,t), v_3(x,t))$
is a velocity perturbation. 

Substituting $u = A\cdot x + v$ in (\ref{eqn:NSE}) and
linearizing in $v$ (\ie omitting $(v \cdot \nabla) v$) yields
\begin{equation}
\begin{split}
\frac{\partial v}{\partial t} + \Bigl((A\cdot x)\cdot \nabla \Bigr) v +
  A\cdot v+ \nabla p  = &\; \nu \nabla^2 v,\\
\nabla \cdot v = &\; 0, \quad
\end{split}
\label{veq}
\end{equation}
or, in indicial notation,
\begin{equation}
\begin{split}
\frac{\partial v_l}{\partial t} +
 A_{mn}x_n \frac{\partial  v_l}{\partial x_m}
  + A_{lm} v_m + \frac{\partial p}{\partial x_l}  =& \;
 \nu \frac{\partial^2 v_l}{\partial x_m  \partial x_m},\\
 \frac{\partial v_l}{\partial x_l} = & \; 0.
\label{veqind}
\end{split}
\end{equation}
We will use the 
convention of summation over repeated indices $l,m,n=(1,2,3)$.
The self-advection of the background field is the gradient of
a potential  that has been absorbed into the pressure, since
 $dA/dt + A\cdot A$ is symmetric.
Incompressibility of the base field requires $\mbox{tr}(A)= A_{ll}=0$.
An equation for the pressure can be derived, in a standard manner,
by taking the divergence of the momentum equation in (\ref{veqind})
and using incompressibility to obtain 
\begin{equation}
\nabla^2 p= \frac{\partial^2 p}{\partial x_l  \partial x_l}
= - 2 A_{ml}  \frac{\partial  v_l}{\partial x_m}.
\label{pressure}
\end{equation}

Linear problems of the form (\ref{veq}) have been considered
by many authors since Kelvin. There has been a renewed
interest in such analyses in more recent times linked 
to a simple model of the elliptical instability (see \cite{RRK02}
and references therein). 

\textbf{Kelvin modes.} 
Kelvin noticed that a generalized Fourier analysis could 
be used to solve equation (\ref{veq}). He proposed 
to look for solutions of the form 
\begin{equation}
v(x,t) = \hat{v}(t) e^{i k(t) \cdot x}
\label{kelvin}
\end{equation}
\ie Fourier modes with time-dependent wave-vectors 
$k \in \Real^3$.
The motivation for this ansatz is that a Fourier mode
initial condition proportional to $e^{i k \cdot x}$
is rotated and stretched uniformly by the background
field with uniform gradients. Hence, it remains in
the form of a Fourier mode, albeit with an evolving 
wavevector $k$. Substituting (\ref{kelvin}) together with $p(x,t)=
\hat{p}(t) e^{i k(t) \cdot x}$ into (\ref{veq}), and using
(\ref{pressure}) to eliminate the pressure, leads to 
the coupled ordinary differential equations:
\begin{gather}
\frac{dk}{dt} = - k \cdot A,
\label{kel1}\\
\frac{d\hat{v}}{dt} = - \nu |k|^2 \hat{v} -A \cdot \hat{v} +
 2 \frac{k}{|k|^2} \left(k\cdot A \cdot v\right). \label{kel2}
\end{gather}
The incompressibility constraint $\nabla\cdot v=0$ requires
$k(t)\cdot \hat{v}(t)=0$. This is satisfied  automatically 
since the pressure is determined from (\ref{pressure}), provided
that  the
initial conditions are such that $k(0)\cdot \hat{v}(0)=0$.
The most interesting solutions of these equations, perhaps, occur
when the base flow \mbox{$A\cdot x$} has closed streamlines. 
Then $k(t)$ is oscillatory
and $\hat{v}(t)$ can grow exponentially through
a  parametric instability \cite{RRK02}.

\textbf{Fourier transform.}
Alternatively, one can proceed with a direct Fourier transform
of the equations. Let $\hat{v}(k,t)$ be the Fourier transform 
of the vector
\begin{equation}
v(x,t) = \int_{\Real^3} \hat{v}(k,t) \, e^{i k \cdot x} dk,
\end{equation}
then 
\begin{equation}
A_{mn} x_n \frac{\partial v_l(x,t)}{\partial x_m}  = 
-\int_{\Real^3} A_{mn} \frac{\partial}{\partial k_n}
\Bigl(k_m \hat{v_l}(k,t)\Bigr) \, e^{i k \cdot x} dk,
\end{equation}
which, in our case, simplifies to 
\begin{equation}
A_{mn} x_n \frac{\partial v_l(x,t)}{\partial x_m}  = 
-\int_{\Real^3} k_m A_{mn} \frac{\partial \hat{v}_l(k,t)}{\partial k_n}
\, e^{i k \cdot x} dk,
\end{equation}
because $A_{mm}=0$.
It follows that the Fourier transform of equation (\ref{veqind}) reads

\begin{equation}
\frac{\partial \hat{v}_l(k,t)}{\partial t} 
- k_m A_{mn} \frac{\partial \hat{v}_l}{\partial k_n}
= -\nu |k|^2 \hat{v}_l - A_{lm} \hat{v}_m +
 2 \frac{k_l}{|k|^2} \bigl(k_m A_{mn} \hat{v}_n\bigr), 
\label{FTveq}
\end{equation}
where the pressure and the incompressibility constraint have
 been eliminated using (\ref{pressure}) provided the 
initial conditions satisfy $k_l \hat{v}_l(k,0)=0$, $\forall k$.
Solving equation (\ref{FTveq}) by the method of characteristics,
we directly recover Kelvin modes and  equations (\ref{kel1}),
(\ref{kel2}).

Equation (\ref{FTveq}) is identical to Dinaburg and Sinai's equation
(11) except for a different definition of the matrix $A$, theirs is
minus the transpose of ours: $A^{(DS)}= -A^T$. 
Note that equation (\ref{FTveq}) has been obtained without imposing
the point symmetry $v(x,t)=-v(-x,t)$, assumed by
Dinaburg and Sinai, although the base flow \mbox{$A\cdot x$} does
satisfy that symmetry. The point symmetry 
implies that $\hat{v}(k,t)=-\hat{v}(-k,t)$. Now,
$\hat{v}(k,t)=\overline{\hat{v}(-k,t)}$ since $v(x,t)$ is real,
 where the overline denotes complex conjugate.
Hence, $\hat{v}(k,t)$ must be pure imaginary if the point symmetry is
imposed.  So Dinaburg and Sinai's $v(k,t)$ is related to our
$\hat{v}(k,t)$ as $v^{(DS)}(k,t) = -i \hat{v}(k,t)$.

 Dinaburg and Sinai define
\begin{equation}
A^{(DS)}_{mn}= -\int i k_m \hat{v}_n(k,t) dk  = 
- \left. \frac{\partial v_n(x,t)}{\partial x_m} \right|_{x=0},
\label{ADS}
\end{equation}
while, if we write our base flow as $U(x,t)=A(t)\cdot x$, then
 $U_m = A_{mn} x_n$ and 
\begin{equation}
A_{mn}= \frac{\partial U_m(x,t)}{\partial x_n}.
\end{equation}
This explains the differences in the definitions of the velocity
gradient $A$.
It also points to a basic inconsistency of the Dinaburg-Sinai model.
Equation (\ref{FTveq}) corresponds to 
the distortion of small scales by an infinitely large scale flow 
with uniform gradient, but Dinaburg and Sinai define the gradient
of the infinitely large scale velocity field as the local gradient
 at $x=0$ resulting from the linear superposition of all the small
scale gradients.
Considering two distinct Fourier modes 
with wavevectors $k^{(1)}$ and $k^{(2)}$, for instance, with 
$ |k^{(1)}| < |k^{(2)}|$, it does not make sense to have the small scale
$k^{(2)}$ participating in the uniform large scale distortion of the 
larger scale $k^{(1)}$. Furthermore, the mode $k^{(1)}$ does not 
self-distort because of the incompressibility constraint 
(so $(v \cdot \nabla) v =0$ for a single Fourier mode).
Hence, only larger scales should contribute to the approximation
of uniform distortion of a given small scale. In other words, 
only $|k'| \ll |k|$
should contribute to the large scale gradient distorting the 
mode with Fourier wavevector $k$. This can be seen also at a technical 
level in the derivation of the Dinaburg-Sinai model. In their 
treatment of the convolution integral representing the
Fourier transform of the Navier-Stokes nonlinearity,
Dinaburg and Sinai make the assumption that
the integral is dominated by highly nonlocal interactions,
\ie by the domains  $|k'| \ll |k|$ and $|k-k'| \ll |k|$.
Considering the domain $|k'| \ll |k|$, they make the following
type of approximation, for instance,
\begin{equation}
 \int_{\Real^3}   k'_n u_m(k') u_n(k-k')   dk' \approx
u_n(k)  \int k'_n  u_m(k') dk'
\label{convol}
\end{equation}
(see eqn.\ (6) in \cite{DS01}).
The domain of integration for the integral on the right-hand side
should be restricted to the ball $B_{\epsilon}=$ 
$\{k' \in \Real^3: |k'| < \epsilon \} $, 
with $\epsilon \ll |k|$, for self-consistency, but Dinaburg and Sinai
do not specify the domain of integration. In their later studies
of finite dimensional approximations, they sum over all modes,
 in other words, they integrate over $\Real^3$ instead of
$B_{\epsilon}$. If the integral, and the finite dimensional sums,
were correctly restricted, the model would not be closed, or
the infinitely large scale flow would have to be specified and the model would
become linear and identical to (\ref{FTveq}), with $A(t)$ specified
independently from the small scales.

It is natural to wonder whether one could replace the Dinaburg-Sinai
model by a mean field theory where the large
scale flow results from the nonlinear interactions of the 
small scales. The base field
$U=A\cdot x$ is a very singular $k=0$ mode whose generalized 
Fourier transform $\hat{U}(k) = i (A \cdot \nabla_k) \delta(k)$.
Here $\delta(k)=\delta(k_1)\delta(k_2)\delta(k_3)$ is a product
of Dirac delta functions and $\nabla_k$ is the gradient operator in $k$-space
with  $\partial\delta(k)/\partial k_1 = \delta'(k_1) \delta(k_2) \delta(k_3)$,
where $\delta'(k_1)$ is the generalized derivative of the delta function,
and similarly for derivatives with respect to $k_2$ and $k_3$.
The only nonlinear interactions that can create a $k=0$ mode consist
of any 
mode $k'$ interacting with its complex conjugate $-k'$, but such interactions
vanish because of incompressibility. For instance, the convolution
integral on the left-hand side of (\ref{convol}) when $k=0$,
\begin{equation}
 \int_{\Real^3}   k'_n u_m(k') u_n(-k')   dk' =0
\label{nozero}
\end{equation}
for any regular $u(k')$ because $k'_n u_n(-k')=0$ from incompressibility.
The vanishing of such interactions is also related to Galilean invariance
since a $k=0$ mode could also correspond to a constant velocity.

\end{document}